%

\input amstex
\documentstyle{amsppt}
\magnification=1200
 \vsize19.5cm
  \hsize13.5cm
   \TagsOnLeft
\pageno=1
\baselineskip=15.0pt
\parskip=3pt

\def\p{\partial}
\def\noo{\noindent}
\def\eps{\varepsilon}

\def\Om{\Omega}

\def\pom{{\p \Om}}
\def\bom{{\overline\Om}}
\def\R{\Bbb R}

\def\det{\text{det}}

\def\ol{\overline}
\def\lan{\langle}
\def\ran{\rangle}
\def\D{\nabla}
\def\phi{\varphi}

\nologo \NoRunningHeads

\topmatter

\title{On strict convexity and $C^1$ regularity\\
 of potential functions in optimal transportation}\endtitle

\author{Neil Trudinger\ \ \ \  Xu-Jia Wang}\endauthor

\affil{Centre for Mathematics and Its Applications\\
The Australian National University }\endaffil

\address Neil S. Trudinger, Centre for Mathematics and Its Applications,
Australian National University,  Canberra, ACT 0200, Australia
\endaddress
\email neil.trudinger\@maths.anu.edu.au  \endemail

\address
Xu-Jia Wang, Centre for Mathematics and its Applications,
Australian National University, Canberra ACT 0200, Australia
\endaddress
\email wang\@maths.anu.edu.au \endemail

\thanks
This work was supported by the Australian Research Council
\endthanks

\abstract {This note concerns the relationship between conditions
on cost functions and domains and the convexity properties of
potentials in optimal transportation and the continuity of the
associated optimal mappings. In particular, we prove that if the
cost function satisfies the condition (A3), introduced in our
previous work with Xinan Ma, the densities and their reciprocals
are bounded and the target domain is convex with respect to the
cost function, then the potential is continuously differentiable
and its dual potential strictly concave with respect to the cost
function. Our result extends, by different and more direct proof,
similar results of Loeper proved by approximation from our earlier
work on global regularity. }
\endabstract


\endtopmatter

\vskip-10pt

\document

\baselineskip=12.5pt
\parskip=3pt

\centerline{\bf \S 1. Introduction}

\vskip10pt

We continue our investigation on the regularity of potential
functions in the optimal transportation problem [MTW, TW]. In this
paper we prove strict convexity and $C^1$ regularity of potential
functions for non-smooth densities. The strict convexity and $C^1$
regularity for solutions to the Monge-Ampere equation were
established by Caffarelli [C1]. For the reflector design problem,
which is a special optimal transportation problem [W2], these
results were obtained in [CGH]. For more general optimal
transportation problems, the $C^{1, \alpha}$ regularity has been
obtained by Loeper [L]. In this paper we prove the strict
convexity for potential functions, and obtain the $C^1$ regularity
under weaker conditions on the domains. We also use our results to
plug a gap in [MTW] pertaining to the use of a comparison argument
in the proof of interior regularity.

Let $\Om, \Om^*$ be two bounded domains in $\R^n$, and $f, g$ be
two nonnegative integrable functions on $\Om, \Om^*$ satisfying
the mass balance condition
$$\int_\Om f=\int_{\Om^*} g. \tag 1.1$$
Let $(u, v)$ be potential functions to the optimal transportation
problem, namely $(u, v)$ is a maximizer of
$$\sup\{I(\phi, \psi):\ (\phi, \psi)\in K\},\tag 1.2$$
where
$$\align
& I(\phi, \psi)
 =\int_\Om f(x)\phi(x)+\int_{\Om^*} g(y)\psi(y), \tag 1.3\\
& K=\{(\phi, \psi)\in C^0(\Om)\times C^0(\Om^*):
                      \ \ \phi(x)+\psi(y)\le c(x, y)\}.\\
\endalign$$
We assume that the cost function is smooth, $c\in
C^\infty(\R^n\times \R^n)$, and satisfies (A1)-(A3) below. It is
known [C4, GM] that there is a maximizer $(u, v)$  to (1.2), which
is also unique up to a constant if $f, g$ are positive.  The
potentials $(u, v)$ are semi-concave and satisfy
$$\align
u(x) &= \inf_{y \in {\Omega}^{*}}\{c(x,y) - v(y)\},\tag 1.4\\
v(y) &= \inf_{x \in \Omega}\{c(x,y) - u(x)\}.\\
\endalign $$
The optimal mapping $T:\ x\in\Om\to y\in\Om^*$ can be determined,
a.e. in $\Om$,  by
$$Du(x)=D_{x}c(x,y),\tag 1.5$$
where $y$ also attains the infimum in (1.4). If $u\in C^2(\Om)$,
it satisfies the equation
$$ \det(D^{2}_{x}c - D^{2}u) = |\det c_{i,j}|
        \frac{f}{g\cdot T}\ \ \text{in}\ \ \Om, \tag 1.6,$$
where $c_i=\p_{x_i} c$, $c_{i, j}=\p_{x_i}\p_{y_j} c$. The main
result of this note is

\proclaim{Theorem 1} Suppose that $\Om^*$ is $c$-convex with
respect to $\Om$, that the cost function $c$ satisfies (A1)-(A3),
and that $f, g$ satisfy
$$C_1\le f, g\le C_2  \tag 1.7$$
for some positive constants $C_1, C_2$. Then $v$ is strictly
$c$-concave and $u$ is $C^1$ smooth.
\endproclaim

We refer the reader to \S2.1 for definitions of various convexity
notions relative to cost functions. By approximation and the
uniqueness of potential functions (when $f, g>0$), condition (1.7)
can be weakened to $0\le \frac{f}{g\cdot T}<C$ and $f, g>0$. From
the $C^1$ smoothness of $u$ it follows that the optimal mapping
$T$ is continuous.

For the reflector design problem (in the far field case), Theorem
1 was obtained in [CGH]. For the optimal transportation problem,
the $C^1$ smoothness of $u$ essentially belongs to Loeper [L]. His
proof uses approximation by globally smooth solutions to the
optimal transportation problem, established in [TW], and
accordingly assumes additional conditions such as the
$c$-convexity of both domains $\Om$ and $\Om^*$. See also remarks
after Corollary 1 in \S 2.8.  Here we give a more direct proof.
Our proof is also completely different from that in [C1], in which
Caffarelli proved the strict convexity and $C^1$ regularity for
solutions to the Monge-Ampere equation with constant boundary
condition. For higher regularity, the interior and global $C^{2,
\alpha}$ estimates for solutions of (1.6), in the case of the
quadratic cost function, were established in [C2, C3, U1], and
earlier in [D] for $n=2$. For cost functions satisfying (A1)-(A3)
below, the regularity of potential functions was obtained in [MTW,
TW].

By approximation and the $C^1$ regularity in Theorem 1, it follows
that the potential function $u$ is {\it fully} $c$-concave if $f,
g>0$ and $c$ satisfies (A1)-(A3). In particular, every local
$c$-support is a global one. See \S 2.8 for more discussion.

The assumption (A1)-(A3) are as follows.

\noo (A1) For any $x, z\in\R^n$, there exists a unique $y=y(x, z)$
such that $D_{x}c(x,y) = z$.

\noo (A2) For any $x, y\in\R^n$, $\det \{c_{i, j}(x, y)\} \ne 0.$

\noo (A3) For any $x\in\Om, y\in\Om^*$, and $\xi, \eta\in\R^n$
with $\xi  \perp \eta$,
$${\sum}_{i, j, k, l, p, q, r, s}
  (c^{p,q}c_{ij,p}c_{q,rs} - c_{ij,rs})c^{r, k}c^{s, l}
  \xi_i\xi_j\eta_k\eta_l\ge c_0|\xi|^2|\eta|^2, \tag 1.8$$
where $c_0$ is a positive constant, and $(c^{i,j})$ is the inverse
matrix of $(c_{i,j})$.

We also assume the above conditions hold after exchanging $x$ and
$y$. Formula (1.8) is equivalent to
$$\p_{z_kz_l} c_{ij}(x, y) \xi_i\xi_j\eta_k\eta_l
                   \le -c_0|\xi|^2|\eta|^2\tag 1.9$$
for all $x\in\Om, y\in\Om^*$ and $\xi\perp \eta\in\R^n$, where
$y=y(x, z)$ is given by (A1), which is smooth in $x$ and $z$ by
(A2). Note that for global regularity in [TW] and the subsequent
application in [L], condition (A3) can be relaxed to its
degenerate form, $c_0=0$, called A3w in [TW].

We divide the proof of Theorem 1 into several short sections. We
first introduce in \S2.1 various notions of convexities and
concavities relative to the cost function $c$. We then indicate in
\S 2.2 a geometric property of (A3) (see also [L]). In \S2.3 we
give an analytic formulation of the $c$-convexity of domains. A
geometric characterization of $c$-convex domains (under condition
A3) is given in \S2.4.  In \S 2.5 we prove that a local
$c$-concave function is fully $c$-concave if the domain is
$c$-convex and the cost function $c$ satisfies (A3). This is a key
ingredient in the proof of Theorem 1. In \S 2.6 we show that $u$
is $C^1$ if and only if $v$ is strictly $c$-concave. We then prove
Theorem 1 in \S 2.7. Various remarks are given in \S 2.8.

\vskip30pt

\centerline{\bf 2. Proof of Theorem 1}

\vskip10pt

{\bf 2.1. Convexities relative to cost functions} [MTW]. Let $u$
be a semi-concave function in $\Om$, namely $u-C|x|^2$ is concave
for a large positive constant $C$. The supergradient $\p^+u$ [GM]
and $c$-supergradient  $\p^+_c u$ are defined by
$$\align
\p^+u(x_0) &= \{p\in\R^n:\ \
          u(x)\le u(x_0)+p\cdot (x-x_0)+o(|x-x_0|)\},\tag 2.1\\
\p^+_cu(x_0) &=\{y\in\R^n:\ \
          u(x)\le c(x, y)-c(x_0, y)+u(x_0)+o(|x-x_0|)\}\tag 2.2\\
          \endalign$$
for $x$ near $x_0$, where $x_0\in\Om$. For a set $E\in \Om$, we
denote $\p^+u(E)=\cup_{x\in E} \p^+u(x)$ and
$\p^+_cu(E)=\cup_{x\in E}\p^+_cu(E)$. By (1.5) we have
$$c_x(x_0, \p^+_cu(x_0))=\p^+u(x_0).\tag 2.3$$
Note that $\p^+u(x_0)$ is a closed, convex set. Hence
$\p^+_cu(x_0)$ is closed and $c$-convex with respect to $x_0$.

We may extend the above mappings to boundary points. Let
$x_0\in\pom$ be a boundary point, we denote
$\p^+u(x_0)=\{p\in\R^n:\ p=\lim_{k\to\infty} p_k\}$, where $p_k\in
\p^+u(x_k)$ and $\{x_k\}$ is a sequence of interior points of
$\Om$ such that $x_k\to x_0$, and let $\p^+_cu(x_0)$ be given by
(2.3).

The $c$-normal mapping $T_u$ is defined by
$$T_u(x_0)=\{y\in\R^n:\ \
  u(x)\le c(x, y)-c(x_0, y)+u(x_0)\ \text{for all}\ x\in\Om\}.\tag 2.4$$
Note that $T_u(x_0)\subset \p^+_c u(x_0)$.

{\it $c$-support}: Let $u$ be a semi-concave function in $\Om$. A
{\it local $c$-support} of $u$ at $x_0\in\bom$ is a function of
the form
$$h=c(\cdot, y_0)+a_0, $$
where $a_0$ is a constant and $y_0\in\R^n$, such that $u(x_0) =
h(x_0)$ and $u(x) \le h(x)$ near $x_0$. If $u(x) \le h(x)$  for
all $x\in \Om$, then $h$ is a {\it global $c$-support} (or
$c$-support for short) of $u$ at $x_0$. If $h$ is a local
$c$-support of $u$ at $x_0$, then $y_0\in \p^+_cu(x_0)$ and
$a_0=u(x_0)-c(x_0, y_0)$.

{\it $c$-concavity of functions}: We say a semi-concave function
$u$ is {\it locally $c$-concave} if for any point $x_0\in \bom$
and any $y\in \p^+_cu(x_0)$, $h=c(\cdot, y)-c(x_0, y)+u(x_0)$ is a
local $c$-support of $u$ at $x_0$. We say $u$ is {\it $c$-concave}
if for any point $x_0\in \bom$, there exists a global $c$-support
at $x_0$ in $\Om$. We say $u$ is {\it strongly $c$-concave} if it
is both locally $c$-concave and $c$-concave. We say $u$ is {\it
fully $c$-concave} if it is locally $c$-concave and every local
$c$-support of $u$ is a global $c$-support. We say $u$ is {\it
strictly $c$-concave} if it is fully $c$-concave and every
$c$-support of $u$ contacts its graph at one point only.

{\it $c$-segment}: A set of points $\ell\subset\R^n$ is a
$c$-segment with respect to a point $y_0\in\R^n$ if $D_yc(\ell,
y_0)$ is a line segment in $\R^n$.

{\it $c$-convexity of domains}: We say a set $U$ is $c$-convex
with respect to another set $V$ if the image $c_y(U, y)$ is convex
for each $y\in V$. Equivalently, $U$ is $c$-convex with respect to
$V$ if for any two points $x_0, x_1\in U$ and any $y\in V$, the
$c$-segment relative to $y$ connecting $x_0$ and $x_1$ lies in
$\Om$. By (2.3), a $c$-convex domain is topologically a ball.

By definition, a $c$-concave function $u$ can be represented as
[GM]
$$u(x)=\inf \{c(x, y)-a(y):\ \ y\in T_u(\Om)\},$$
where $a$ is a function of $y$ only.  By (1.4), a potential
function $u$ is $c$-concave [GM]. Our Theorem 1 implies that $u$
is furthermore fully $c$-concave under assumption (A3). If $u$ is
$C^1$, then local $c$-support is unique and $c$-concavity is
equivalent to full $c$-concavity. We also remark that a potential
function may not be locally $c$-concave in general.

Similarly we can define $c^*$-segment, $c^*$-support,
$c^*$-convexity and $c^*$-concavity by exchanging variables $x$
and $y$ [MTW]. In this paper, we will generally omit the
superscript $*$ when the meaning is clear.

\vskip10pt

{\bf 2.2. A geometric property of (A3)}. Let $y_0, y_1$ be two
points in $\ol{\Om^*}$. Let $\widehat{y_0y}_1$ be the $c$-segment
relative to a point $x_0\in\Om$, connecting $y_0$ and $y_1$. By
definition,
$$\widehat{y_0y}_1
   =\{y_t:\  c_x(x_0, y_t)=p_t,\ t\in [0, 1]\},\tag 2.5$$
where $p_t=tp_1+(1-t)p_0$, and $p_0=c_x(x_0, y_0)$, $p_1=c_x(x_0,
y_1)$.  Let
$$\align
h_i(x)& =c(x, y_i)-a_i,\ \ \ i=0, 1, \tag 2.6\\
h_t(x)& =c(x, y_t)-a_t,\\
\endalign $$
where $t\in (0, 1)$, $a_0, a_1$ and $a_t$ are constants such that
$h_0(x_0)=h_1(x_0)=h_t(x_0)$. Suppose (A3) holds. Then for $x\ne
x_0$, near $x_0$, we have the inequality
$$h_t(x)> \min\{h_0(x), h_1(x)\},\tag 2.7$$
which is crucial for the remaining analysis of this paper.

Inequality (2.7) follows from (1.9). Indeed, by a rotation of
axes, we assume that $p_1=p_0+\delta e_n$, where $\delta$ is a
positive constant and $e_n=(0, \cdots, 0, 1)$ is the unit vector
in the $x_n$-axis. By (1.9),
$$\frac{d^2}{dt^2} c_{ij}(x, y(x,p_t))\xi_i\xi_j
                       \le -c_0\delta^2 \tag 2.8$$
for any unit vector $\xi$ orthogonal to $e_n$, where $y=y(x_0,
p_t)$ is given in (A1). Now (2.7) follows from (2.8); for details
see [L].

\vskip10pt

{\bf 2.3. An analytic formulation of the $c$-convexity of
domains}. If $\Om$ is $c$-convex with respect to $\Om^*$, by
definition, $c_y(\Om, y)$ is convex for any $y\in\Om^*$. Suppose
$0\in\pom$ and locally $\pom$ is given by
$$x_n=\rho(x')\tag 2.9$$
with $D\rho(0)=0$ such that $e_n$ is the inner normal at $0$,
where $x'=(x_1, \cdots, x_{n-1})$. Then at $(0, y)$ (for a fixed
point $y$),
$$\align
\p_{x_i} c_y & =c_{i, y}+c_{n, y}\rho_{x_i},\\
\p_{x_ix_j} \lan c_y, \gamma\ran
 & =\lan c_{ij, y}, \gamma\ran
  +\lan c_{n, y}, \gamma\ran\rho_{x_ix_j}\ge 0,\tag 2.10\\
\endalign $$
where $\gamma$ is the inner normal of $c_y(\Om, y)$ at $c_y(0,
y)$. We may write (2.10) explicitly,
$$c_{ij, y_l}\gamma_l+c_{n, y_l}\gamma_l \rho_{ij}\ge 0.$$
Make the linear transformation
$$\hat y_k=a_{kl}y_l.$$
Then we have
$$c_{ij, \hat y_k}a_{kl} a^{-1}_{ml}\hat \gamma_m
  +c_{n, \hat y_k}a_{kl} a^{-1}_{ml}\hat \gamma_m \rho_{ij}\ge 0.\tag 2.11$$
Let $a_{ij}= c_{i, j}(0, y)$. Then $c_{x_i, \hat y_j}=\delta_{ij}$
and $\hat \gamma=e_n$. We obtain
$$c_{ij, \hat y_n}+\rho_{x_ix_j}\ge 0,\tag 2.12$$
which is equivalent to
$$c_{ij, y_k}c^{k, n}+\rho_{x_ix_j}\ge 0.\tag 2.13$$

Let $\phi\in C^2(\bom)$ be a defining function of $\Om$. That is
$\phi=0$, $|\D \phi|\ne 0$ on $\pom$ and $\phi<0$ in $\Om$. From
(2.13), we obtain an analytic formulation of the $c$-convexity of
$\Om$ relative to $\Om^*$ [TW],
$$[\phi_{ij}(x)-c^{k, l}(x, y)c_{ij, k}(x, y)\phi_l(x)]\ge 0
 \ \ \ \forall\ x\in\pom, y\in\Om^* .\tag 2.14$$
Conversely, if $\Om$ is simply connected and (2.14) holds, then
$\Om$ is $c$-convex. Following [TW], we call $\Om$ uniformly
$c$-convex with respect to $\Om^*$ if the matrix in (2.14) is
uniformly positive.

\vskip10pt

{\bf 2.4.  Geometric properties of the $c$-convexity of domains}.
Let $y_0, y_1$ be any two given points in $\ol{\Om^*}$. Denote
$$\Cal N=\Cal N_{y_0, y_1, a}
        =\{x\in\R^n:\ c(x, y_0)=c(x, y_1)+a\},\tag 2.15$$
where $a$ is a constant. Assume that the origin $0\in\Cal N$ and
locally $\Cal N$ is represented as
$$x_n=\eta(x')\tag 2.16$$
such that $\eta(0)=0$ and $D\eta(0)=0$ (obviously $\eta$ also
depends on $y_0, y_1$ and $a$). Then
$$c(x', \eta(x'), y_0)=c(x', \eta(x'), y_1)+a. \tag 2.17$$
Differentiating (2.17) gives
$$\align
c_i(0, y_0)+c_n\eta_i &=c_i(0, y_1)+c_n\eta_i,\\
c_{ij}(0, y_0)+c_n\eta_{ij} &=c_{ij}(0, y_1)+c_n\eta_{ij}.\tag 2.18\\
\endalign $$
We obtain
$$[c_n(0, y_1)-c_n(0, y_0)]\eta_{ij}
      +[c_{ij}(0, y_1)-c_{ij}(0, y_0)]=0. \tag 2.19$$

Now let $y_0$ be fixed but let $y_1$ and  $a$ vary in such a way
that $y_1\to y_0$ and the set $\Cal N=\Cal N_{y_0, y_1, a}$, given
by (2.15), is tangential to $\{x_n=0\}$, namely $\eta(0)=0$ and
$D\eta(0)=0$. By a linear transform as in \S 2.3 we assume that
$c_{x_i, y_j}=\delta_{ij}$ at $x=0$ and $y=y_0$. Then
$$\frac {y_1-y_0}{|y_1-y_0|}\to e_n. $$
We obtain
$$c_{n, n}(0, y_0)\eta^0_{ij}+c_{ij, n}(0, y_0)=0,  \tag 2.20$$
where $\eta^0$ is the limit of $\eta_{y_0, y_1, a}$ as $y_1$ and
$a$ vary as above.

Now let $\Om$ be $c$-convex, (uniformly $c$-convex), with respect
to $\Om^*$. Suppose that $\pom$ is given by (2.9) with
$D\rho(0)=0$ so that $\pom$ is tangential to $\Cal N$ at the
origin.  Then by (2.12)and (2.20) we obtain
$$D^2 (\rho-\eta^0)\ge 0,\ \ (>0), \ \ \text{at}\ 0. \tag 2.21$$
From (2.21) we obtain some useful geometric properties of
$c$-convex domains, assuming that $c$ satisfies (A3).

First, if $\Om$ is $c$-convex with respect to $\Om^*$, then for
any compact subset $G\subset\Om^*$, $\Om$ is uniformly $c$-convex
with respect to $G$. Indeed, let $y_0, y_1$ be two points in
$\ol{\Om^*}$. Let $p_0=c_x(0, y_0)$, $p_1=c_x(0, y_1)$.  For $t\in
(0, 1)$, let $p_t=tp_1+(1-t)p_0$ and $y_t$ satisfy $c_x(0,
y_t)=p_t$. Then the set $\{p_t:\ 0\le t\le 1\}$ is a line segment
and the set $\{y_t:\ 0\le t\le 1\}$ is a $c$-segment. Suppose
$p_1=p_0+\delta e_n$ for some $\delta>0$. Let $h_t(x)=c(x,
y_t)+a_t$, where $a_t=c(0, y_0)-c(0, y_t)$ such that
$h_t(0)=h_0(0)$ for all $t\in [0, 1]$. Let $\Cal N_t=\{x\in\R^n:\
h_t(x)=h_0(x)\}$ and suppose that near $0$, $\Cal N_t$ is given by
$x_n=\eta_t(x')$. Since $c_x(0, y_t)=p_t$, we see that $\Cal N_t$
is tangential to $\{x_n=0\}$, namely $D(\eta_t-\eta_{t'})(0)=0$.
By (2.8) we have furthermore the monotonicity formula
$$D^2(\eta_t-\eta_{t'})(0) >0 \tag 2.22$$
for any $t>t'$ and $t, t'\in [0, 1]$. Geometrically it implies
that $\Cal N_t$ lies above $\Cal N_{t'}$ if $t>t'$, namely
$\eta_t(x')\ge \eta_{t'}(x')$ for $x'$ near $0$, and equality
holds only at $x'=0$. Consequently we obtain from (2.21) the
strict inequality
$$D^2 (\rho-\eta)>0\ \ \ \text{at}\ 0. \tag 2.23$$
From (2.23) we obtain the above mentioned property.

Next, for any $y_0, y_1\in\overline{\Om^*}$, if $\Cal N$ (given in
(2.15)) is tangent to $\pom$ at some point $x_0$ and if $\Om$ is
$c$-convex with respect to $\Om^*$, then the whole domain $\Om$
lies on one side of $\Cal N$. Indeed, we may assume $x_0=0$,
locally $\pom$ and $\Cal N$ are given respectively by (2.9) and
(2.16), such that $D\rho (0)=D\eta(0)=0$. By (2.23),
$\rho(x)>\eta(x)$ for $x$ near $0$, $x\ne 0$. Denote
$$U=U_{y_0, y_1, a}=\{x\in\R^n:\ c(x, y_0)< c(x, y_1)+a\}$$
so that $\Cal N=\p U$. If $\Om$ does not lie on one side of $\Cal
N$, namely $\Om$ is not contained in $U$, then $\bom-U$ contains
two disconnected components (one is the origin). Since $\pom$ is a
closed, compact hypersurface, we decrease the constant $a$
(shrinking the set $U$) until a moment when two components of
$\bom-U$ meet each other at some point $x^*\in\pom$. But since
$\Cal N$ is tangent to $\pom$ at $x^*$, we reach a contradiction
by (2.23).

It follows that if $\Om$ is $c$-convex with respect to $\Om^*$,
then
$$\Om=\bigcap U_{y_0, y_1, a} , \tag 2.24$$
where the intersection is for all $y_0, y_1\in\ol{\Om^*}$ and
constant $a$ such that $U_{y_0, y_1, a}\supset \Om$. However, we
don't know if the converse is true, namely whether $\Om$ is
$c$-convex with respect to $\Om^*$ if it is given by (2.24).

The above properties also extend to cost functions satisfying A3w
and uniformly $c$-convex domains.

\vskip10pt

{\bf 2.5. Local $c$-support is global}. Let $u$ be a locally
$c$-concave function in $\Om$ with $\p^+_cu(\Om)\subset\Om^*$.
Suppose $\Om$ is $c$-convex with respect to $\Om^*$. Let
$h(x)=h_a(x)=c(x, y_0)+a$ be a local $c$-support of $u$ at $x_0$.
Then $h$ is a global $c$-support of $u$, namely
$$u(x)\le h(x)\ \ \forall\ x\in\Om. \tag 2.25$$

Indeed, if this is not true, then for $\eps>0$ small, the set
$\{x\in\Omega:\ u(x)>h_{a-\eps}(x)\}$ contains at least two
disconnected components. We increase $\eps$ (moving the graph of
$h$ vertically downwards) until at a moment $\eps=\eps_0>0$, two
components first time touch each other at some point $x_0$. If
$x_0$ is an interior point of $\Omega$, by definition
$h_{a-\eps_0}$ cannot be a local $c$-support at $x_0$, which
implies that $y_0\not\in \p^+_cu(x_0)$. We claim that for a
sufficiently small $r>0$, $y_0$ does not lie in the set
$\p^+_cu(B_r(x_0))$ either. Indeed, if $h_k=c(\cdot, y_0)+a_k$ for
$k=1, 2, \cdots$ is a sequence of local $c$-support of $u$ at
$x_k$ and if $x_k, a_k\to x_0, a_0$,  we have $y_0\in
\p^+_cu(x_0)$ as $u$ is semi-concave. Hence $h_0=c(\cdot,
y_0)+a_0$ is a local $c$-support of $u$ at $x_0$. This is a
contradiction.

Therefore $h_{a-\eps_0}$ and $u$ are transversal near $x_0$, and
for $\eps<\eps_0$, close to $\eps_0$, locally the set
$\{x\in\Omega:\ u(x)>h_{a-\eps}(x)\}$ cannot contain two
disconnected components. Hence $x_0$ must be a boundary point of
$\Om$.

In case $x_0$ is a boundary point of $\Omega$, we will also reach
a contradiction by (2.23). Similarly as above, $h_{a-\eps_0}$
cannot be a local $c$-support at $x_0$, namely $y_0\not\in
\p^+_cu(x_0)$. Without loss of generality let us assume that
$x_0=0$ and locally $\pom$ is tangent to $\{x_n=0\}$ such that
$e_n$ is an inner normal of $\Om$ at $0$. As before we also assume
that $c_{x_i, y_j}=\delta_{ij}$ at $x=0$ and $y=y_0$. Let
$p_0=c_x(0, y_0)$. By subtracting a linear function of $x$ from
both $c$ and $u$, we assume that $\{x_n=0\}$ is a tangent plane of
$h$ at $0$. Then we have $p_0=0$,  $u(0)=0$, and $u(x)\le o(|x|)$
as $x\to 0$. Since $y_0\not\in \p^+_cu(x_0)$, we have
$$\beta =\lim_{t\to 0}\frac 1t u(te_n)<0.\tag 2.26$$
Let $p_1=\beta e_n$, $p_t=tp_1+(1-t)p_0$ for $t\in [0, 1]$,   and
let $y_t\in\R^n$ be determined by $c_x(0, y_t)=p_t$. Then $p_1\in
\p^+u(0)$ and $y_1\in \p^+_cu(0)$. Since $u$ is locally
$c$-concave, $c(x, y_1)+b$ is a local $c$-support of $u$ at $0$,
where $b$ is a constant such that $c(0, y_1)+b=u(0)=0$. Denote
$$U=\{x\in\R^n:\ h_{a-\eps_0}(x)>c(x, y_1)+b\}$$
and $\Cal N=\p U=\{x\in\R^n:\ h_{a-\eps_0}(x)=c(x, y_1)+b\}$, such
that $0\in\Cal N$. Since $\Om$ is $c$-convex, by (2.23) we see
that $\Om\subset U$. But since $c(x, y_1)+b$ is a local
$c$-support of $u$ at $0$, we have
$$h_{a-\eps_0}(x)>c(x, y_1)+b\ge u(x)\tag 2.27$$
for $x\in \Om$, near the origin. We reach a contradiction by our
choice of $\eps_0$. This completes the proof of (2.25).

From the above proof, we see that if $h$ is a $c$-support of $u$,
then the contact set $\{x\in\Om:\ h(x)=u(x)\}$ cannot contain two
disconnected components (or points). In other words, the contact
set is connected.

We also remark that if $u$ is a potential function to the optimal
transportation with positive mass distributions $f$ and $g$ so
that $u$ is uniquely determined up to a constant, the boundary
point case $x_0\in\pom$ can be reduced to the case $x_0\in\Om$ by
extending $u$ to larger domains, and it is not necessary to define
the mapping $\p^+u$ and $\p^+_cu$ on boundary points.

Taking account of our remark at the end of the previous section,
we also see that (2.25) extends to A3w costs and uniformly
$c$-convex domains. Consequently we obtain an alternate proof of
the $c$-convexity of the solutions in Section 6 of [TW].
Futhermore by domain approximation, we may then extend (2.25)
further to A3w costs and $c$-convex domains.

\vskip10pt

{\bf 2.6. Potential functions}. Let $(u, v)$ be potential
functions to the optimal transportation problem (1.2). Then for
any point $x_0\in\bom$, by (1.4), there exists a point
$y_0\in\ol{\Om^*}$ such that $u(x_0)+v(y_0)=c(x_0, y_0)$. Hence
$$h(x)=c(x, y_0)-v(y_0)\tag 2.28$$
is a $c$-support of $u$ at $x_0$, and
$$h^*(y)=c(x_0, y)-u(x_0)\tag 2.29$$
is a $c$-support of $v$ at $y_0$.

If $u$ is $C^1$ at $x_0$, it has a unique $c$-support at $x_0$. As
a potential function is semi-concave, it is twice differential
a.e. in $\Om$. Hence $u$ has a unique  $c$-support almost
everywhere.

Next we consider the case when $u$ is not $C^1$ at $x_0$. Let us
first introduce the terminology {\it extreme point}. Let $E$ be a
convex set in $\R^k$. We say a point $p\in \ol E$ is an extreme
point of $E$ if there exists a plane $P$ such that $P\cap\ol E$
contains only the point $p$. It is easy to show, by induction on
dimensions, that any interior point in a convex set can be
expressed as a linear combination of extreme points.

If $u$ is not $C^1$ at an interior point $x_0$, since $u$ is
semi-concave, $\p^+u(x_0)$ is a convex set of dimension $k$ for
some integer $1\le k\le n$. Let $N^e_u(x_0)$ denote the set of
extreme points of $\p^+u(x_0)$. Let
$$T^e_u(x_0)=\{y\in\R^n:\ c_x(x_0, y)\in N^e_u(x_0)\}.\tag 2.30$$
Then for any $y_0\in T^e_u(x_0)$, the function
$$h(x)=c(x, y_0)+a$$
is a global $c$-support of $u$ at $x_0$, where $a$ is a constant
such that $h(x_0)=u(x_0)$.

This assertion follows from a similar one for concave functions,
which can be proved by blowing up the graph of $u$ to a concave
cone. That is for any $p\in T^e_u(x_0)$, there exists a sequence
of $C^1$-smooth points $\{x_k\}$ of $u$, $x_k\to x_0$, such that
the $c$-support of $u$ at $x_k$ converges to a $c$-support of $u$
at $x_0$. Recall that at $C^1$ points, $u$ has a unique global
$c$-support.

From the above assertion and (2.7) it follows that a potential
function is a local $c$-concave function if (A3) is satisfied.
That is for any $y'\in \p^+_cu(x_0)$, the function
$$h(x)=c(x, y')-c(x_0, y')+u(x_0)$$
is a local $c$-support of $u$ at $x_0$. By (1.4), $h$ is a global
$c$-support if
$$u(x_0)+v(y')=c(x_0, y').\tag 2.31$$
We remark that in general a potential function may fail to be
locally $c$-concave (with respect to the definition in \S 2.1), if
(A3w) is violated.

\vskip10pt

From the above assertion it also follows that if $u$ is not $C^1$
at $x_0$, then the function $h^*$ in (2.29) is a $c$-support of
$v$ at any point $y_0\in T^e_u(x_0)$. In other words, $T^e_u(x_0)$
is contained in the contact set
$$\Cal C=\{y\in\Om^*:\ h^*(y)=v(y)\}.\tag 2.32$$
By (2.31), a local $c$-support of $u$ at $x_0$ is a global one if
and only if the contact set $\Cal C$ is $c$-convex (with respect
to $x_0$).

Since $\Om^*$ is $c$-convex with respect to $\Om$, from the
argument in \S 2.5, any local $c$-support of $v$ is a global one.
Hence  $v$ is fully $c$-concave in $\Om^*$. Furthermore, $u$ is
$C^1$ smooth if and only if its dual function $v$ is strictly
$c$-concave.

\vskip10pt

{\bf 2.7 Proof of Theorem 1}. In the proof of Theorem 1 we will
use the Perron lifting. Consider the Dirichlet problem
$$\align
 \det(D^{2}_{y}c - D^{2}w) &= \psi(x, Dw)
    \ \ \ \text{in}\ \ B_r(y_0),\tag 2.33\\
     w & = \phi\ \ \ \text{on}\ \ \p B_r(y_0), \\
 \endalign $$
where $\psi>0$, $\phi, \psi\in C^\infty$. Suppose $r>0$ is
sufficiently small and there is a supersolution to (2.33) (note
that when $r$ is sufficiently small, there is always a
supersolution). From [MTW], there is a solution $w\in C^\infty(\ol
B_r(y_0))$ to (2.33) such that the matrix $(D^{2}_{y}c -
D^{2}w)>0$. By approximation and the interior a priori estimates
[MTW],  there is a solution $w\in C^\infty(B_r(y_0))\cap
C^0(B_r(y_0))$ if $\phi\in C^0$ and $\psi\in C^{1,1}$.  Obviously
$w$ is locally $c$-concave. By \S2.5, $w$ is $c$-concave in $B_r$.

Let $(u,v)$ be the potential functions to (1.2). To prove Theorem
1, it suffices to prove that $v$ is strictly $c$-concave. By
approximation we assume that $f, g$ are positive and smooth. By \S
2.4 we may also assume that $\Om^*$ is uniformly $c$-convex with
respect to $\Om$. From \S2.6, $v$ is fully $c$-concave and every
local $c$-support of $v$ is a global one. Suppose to the contrary
that $v$ is not strictly $c$-concave. Then there is a $c$-support
$h^*$ of $v$ at some point $y_0\in\Om^*$ such that the contact set
$$\Cal C^*=\{y\in\Om^*:\ h^*(y)=v(y)\}$$
contains more than one point.

From the argument in \S2.5, $\Cal C^*$ cannot contain more than
one disconnected component. In other words, $\Cal C^*$ is
connected. Hence for any $r>0$ small, the intersection $\p
B_{2r}(y_0)\cap\Cal C^*$ is not empty. Let $\{B_j: j=1, \cdots,
k\}$ be finitely many balls with radius $r/2$, centered on $\p
B_r(y_0)$, such that $\bigcup_{j=1}^k B_j \supset \p B_r(y_0)$.

Denote $v_0=v$. For $j=1, \cdots, k$, let $v_j$ be the solution of
$$\align
 \det(D^{2}_{y}c - D^{2}w) &= \frac 12 \delta_0
                   \ \ \ \text{in}\ \ B_1,\tag 2.34\\
  w & = v_{j-1}\ \ \ \text{on}\ \ \p B_r(0), \\
 \endalign $$
where
$$\delta_0=\inf |\det c_{i,j}|\frac{f(x)}{g(T(x))},\tag 2.35$$
which is positive by assumption of Theorem 1. Extend $v_j$ to the
whole $\Om^*$ such that $v_j=v_{j-1}$ in $\Om^*-B_j$. Then $\ol
v=:v_k$ is locally $c$-concave in $\Om^*$. From \S2.5, it is fully
$c$-concave.

By the a priori estimates in [MTW], $v_j$ is smooth in $B_j$, for
$j=1, \cdots, k$. By the comparison principle, we have $v_1<v_0$
in $B_1$ and by induction, $v_j<v_0$ in $B_j$ for all $j=2,
\cdots, k$. It follows that $\ol v<v$ near $\p B_r(y_0)$.

We have therefore obtained another fully $c$-concave function $\ol
v$ which satisfies
$$\align
\ol v & =v\ \ \text{in}\ B_{r/2}(y_0)\cup\{\Om^*-B_{3r/2}(y_0)\},\\
\ol v & <v\ \ \text{near}\ \p B_r(y_0).\\
\endalign $$
Hence the contact set $\{y\in\Om^*:\ \ h^*(y)=\ol v(y)\}$ cannot
be connected. But this is impossible from the argument in \S 2.5
(as remarked at the end of \S 2.5). Hence $v$ is strictly
$c$-concave. This completes the proof of Theorem 1.

\vskip10pt

{\bf 2.8. Remarks}. First we have the following result which
follows from the $C^1$ regularity in Theorem 1 and approximation.

\proclaim{Corollary 1} Suppose that $f, g>0$, $f, g\in L^1$, $c$
satisfies (A1)-(A3), and $\Om^*$ is $c$-convex with respect to
$\Om$. Then the potential function $u$ is fully $c$-concave.
\endproclaim

Corollary 1 is a complement to the paper [MTW]. In [MTW] we
introduced a notion of generalized solution to the boundary value
problem (1.6) and proved  interior regularity under the conditions
in Corollary 1 and assuming also the smoothness of $f$ and $g$. If
the potential function $u$ is not fully $c$-concave, a local
$c$-support of $u$ may not be a global one. In such case, the
definition of generalized solution in [MTW] is not proper and the
comparison principle may not hold in arbitrary sub-domains.
Corollary 1 rules out the possibility provided the cost function
satisfies (A1)-(A3) and $\Om^*$ is $c$-convex with respect to
$\Om$, as assumed in [MTW]. However to get the full $c$-concavity
of $u$ we have to prove Theorem 1 first. Clearly a short and
direct proof is desired. But in dimension 2, the full
$c$-concavity is a direct consequence of (2.7).

\proclaim{Corollary 2} Suppose the cost function $c$ satisfies
(A1)-(A3).  Then any potential functions $(u, v)$ defined in the
whole $\R^2$ are fully $c$-concave.
\endproclaim

Indeed, for any point $x_0\in\Om$, if $u$ is $C^1$ at $x_0$, there
is a unique global $c$-support of $u$ at $x_0$. Otherwise, let
$y_0, y_1$ be any two points in $\p^+_cu(x_0)$. Let $h_0(x)=c(x,
y_0)-a_0$ and $h_1(x)=c(x, y_1)-a_1$ be two $c$-supports of $u$ at
$x_0$. Denote $\Cal N=\{x\in\R^n:\ h_0(x)=h_1(x)\}$. $\Cal N$ is a
curve which divides $\R^2$ into two parts, and both are
non-compact. It suffices to show that
$$h_t\ge h_0  \ \ \ (=h_1)\ \ \text{on}\ \ \Cal N,\tag 2.36$$
where $t\in (0, 1)$ and $h_t$ is as in (2.7). But if (2.36) is not
true, by moving the graph of $h_t$ downwards, we see that there is
a constant $b$ and a point $x'\in\Cal N$ such that
$$\align
h_t(x')-b & =  \min\{h_0(x'), h_1(x')\},\tag 2.37\\
h_t(x)-b &\le \min\{h_0(x), h_1(x)\}
     \ \  x\in\Cal N\ \ \text{near}\ x'.\tag 2.38\\
\endalign$$
But this is in contradiction with (2.7) at $x'$. Hence Corollary 2
holds.

Note that for potential functions on bounded domains, by the
uniqueness of potential functions when restricted to $\{f>0\}$, we
see that $u$ (and similarly $v$) is fully $c$-concave when
restricted to $\{f>0\}$.

The proof of Corollary 2 does not extend to higher dimensions, as
we don't know if there is a point $x'\in\Cal N$ such that (2.37)
and (2.38) hold. But Corollary 2 holds on compact manifolds of any
dimension, as the set $\Cal N$ is compact. In particular it holds
for the reflector design problem (in the far field case) [W1, W2].
But for the reflector design problem, a $c$-support is a
paraboloid with focus at the origin and one can also verify (2.36)
directly [CGH].

We remark that if the cost function $c$ does not satisfy (A3),
then (2.7) may not hold. Loeper [L] shows that if condition A3w is
violated, there are potential functions which are not fully
$c$-concave. Furthermore, the potential function $u$ may not be
$C^1$ smooth even if both $f$ and $g$ are positive and smooth.
Note that Loeper's potential functions are the negative of those
here so that our $c$-concavity is equivalent to his $c$-convexity.

\vskip10pt

\noo{\it Remark}. After this paper was finished, we learned that
Kim and McCann found a direct proof of Corollary 1 above. They
proved that under A3 (A3w, resp.), $\frac {d^2}{dt^2}h_t > 0$
($\ge 0$, resp.), from which it follows that the contact set of
the potential function $u$ with its $c$-support is connected,
where $h_t$ was given in (2.6).

 \vskip30pt

\baselineskip=12pt
\parskip=1pt

\Refs\widestnumber\key{ABC}

\item {[C1]} Caffarelli, L.A.,
             A localization property of viscosity solutions to
             the Monge-Amp\`ere equation and their strict convexity,
             Ann. Math., 131(1990), 129-134.

\item {[C2]}  Caffarelli, L.A.,
              The regularity of mappings with a convex potential,
              J. Amer. Math. Soc., 5(1992), 99-104.

\item {[C3]}  Caffarelli, L.A.,
              Boundary regularity of maps with convex potentials II.
              Ann. of Math. 144 (1996), no. 3, 453--496.

\item {[C4]}  Caffarelli, L.A.,
              Allocation maps with general cost functions,
              in {\it Partial Differential Equations and Applications},
              (P. Marcellini, G. Talenti, and E. Vesintini eds),
              Lecture Notes in Pure and Appl. Math., 177(1996),
              pp. 29-35.

\item {[CGH]} Caffarelli, L.A.,  Gutierrez, C. and Huang, Q.,
              On the Regularity of  Reflector Antennas,
              Annals of Math,  to appear.

\item {[D]}   Delano\"e, Ph.,
              Classical solvability in dimension two of the second
              boundary value problem associated with
              the Monge-Amp\`ere operator,
              Ann. Inst. Henri Poincar\'e, Analyse Non Lin\'eaire,
              8(1991), 443-457.

\item {[GM1]}  Gangbo, W.,  McCann, R.J.,
              Optimal maps in Monge's mass transport problem,
               C.R. Acad. Sci. Paris, Series I,  Math.
               321(1995), 1653-1658.

\item {[GM2]}  Gangbo, W.,  McCann, R.J.,
              The geometry of optimal transportation,
              Acta Math.,  177(1996), 113-161.

\item {[GT]}  Gilbarg, D., Trudinger, N.S.,
              Elliptic partial differential equations of second order,
              Springer, 1983.

\item {[L]}  Gregoire Loeper,
             Continuity of maps solutions of optimal transportation problems,
             preprint.

\item {[MTW]} Ma, X.N., Trudinger, N.S., and Wang, X-J.,
              Regularity of potential functions of the optimal
              transportation problem,
              Arch. Rat. Mech. Anal., 177(2005), 151-183.

\item {[TW]}  Trudinger, N.S., and Wang, X-J.,
              On the second boundary value problem for Monge-Amp\`ere type
              equations and optimal transportation, preprint.

\item {[U1]}  Urbas, J.,
              On the second boundary value problem for equations of
              Monge-Amp\`ere type,
              J. Reine Angew. Math., 487(1997), 115-124.

\item {[U2]} Urbas, J.,
              Mass transfer problems,
              Lecture Notes, Univ. of Bonn, 1998.

\item {[V]}  Villani, C.,
             Topics in optimal transportation problem,
             Amer. Math. Soc., 2003.

\item {[W1]}  Wang, X.J.,
              On the design of a reflector antenna,
              Inverse Problems, 12(1996), 351-375.

\item {[W2]}  Wang, X.J.,
             On the design of a reflector antenna II,
             Calc. Var. PDE,  20(2004), 329-341.

\endRefs

\enddocument

\end